\documentclass[reqno,12pt]{amsart}
\usepackage[T1]{fontenc}
\usepackage[utf8]{inputenc}
\usepackage[english]{babel}
\usepackage{amssymb,amsmath,amsthm,amsfonts,xcolor,enumerate,hyperref,comment,longtable,cleveref,mathtools}

\usepackage{times}
\usepackage{cite}
\usepackage{pdflscape}
\usepackage[mathcal]{euscript}
\usepackage{tikz}
\usepackage{cancel}
\usepackage{stmaryrd}
\usepackage{longtable}

\usepackage[a4paper,top=3cm, bottom=3cm, left=3cm, right=3cm]{geometry}

\makeatletter
\@namedef{subjclassname@2020}{%
  \textup{2020} Mathematics Subject Classification}
\makeatother

\theoremstyle{plain}
\newtheorem{theorem}[subsection]{Theorem}

\newtheorem{lemma}[subsection]{Lemma}	
\theoremstyle{definition}
\newtheorem{definition}[subsection]{Definition}
\theoremstyle{remark}

\begin{document}
	
	\title[Local super-derivations of the super Schr\"{o}dinger algebras]{Local super-derivations of the super Schr\"{o}dinger algebras}
	
\author{Alauadinov A.K.,\ Yusupov B.B.}
\address[Amir Alauadinov]{
Ch.Abdirov 1, Department of Mathematics, Karakalpak State University, Nukus 230113, Uzbekistan
\newline
and\newline	
V.I.Romanovskiy Institute of Mathematics, Uzbekistan Academy of Sciences, Univesity Street, 9, Olmazor district, Tashkent, 100174, Uzbekistan		
}
	\email{amir\_t85@mail.ru}

	\address[Bakhtiyor Yusupov]{	
		V.I.Romanovskiy Institute of Mathematics, Uzbekistan Academy of Sciences, Univesity Street, 9, Olmazor district, Tashkent, 100174, Uzbekistan
\newline
and\newline
Department of Physics and Mathematics, Urgench State University, H. Alimdjan street, 14, Urgench
220100, Uzbekistan
}
\email{baxtiyor\_yusupov\_93@mail.ru\ b.yusupov@mathinst.uz}

	
	\begin{abstract} The present paper is devoted to study local super-derivations of the super Schr\"{o}dinger algebras. We prove that every local super-derivation on the super Schr\"{o}dinger algebra is a super-derivation.	
	\end{abstract}


	\maketitle
	
\section{Introduction}

	This notion was  introduced and investigated independently by R.V. Kadison~\cite{Kad} and D.R.
	Larson and A.R. Sourour~\cite{Lar}.
	The above papers gave rise to a series of works devoted to the description of mappings which are close to automorphisms and derivations of $C^\ast$-algebras and  operator algebras.  R.V. Kadison  set out a program of study for local maps in \cite{Kad}, suggesting that local derivations could prove useful in building derivations with particular properties. R.V.Kadison proved in \cite[Theorem A]{Kad} that each continuous local derivation of a von Neumann algebra $M$  into a dual Banach  $M$-bimodule is a derivation. This theorem gave way to  studies on derivations on  $C^\ast$-algebras, culminating with a result due to B.E. Johnson, which asserts that every  local derivation of a  $C^\ast$-algebra  $A$ into a Banach  $A$-bimodule is automatically continuous, and hence is a derivation \cite[Theorem 5.3]{Jon}.

	Let us present a list of finite or infinite dimensional  algebras for which all local derivations are derivations:
	\begin{itemize}
		\item[-] $C^\ast$-algebras, in particular,   the algebra $M_n(\mathbb{C})$ of all square matrices  of order $n$ over the field of complex numbers \cite{Kad, Jon};
		\item[-] the complex polynomial algebra $\mathbb{C}[x]$  \cite{Kad};
		\item[-] finite dimensional simple Lie algebras over an algebraically closed field of characteristic zero \cite{AK};
		\item[-] Borel subalgebras of finite-dimensional simple Lie algebras \cite{YC20};
		\item[-] infinite dimensional Witt  algebras over an algebraically closed field of characteristic zero \cite{CZZ};
		\item[-] Witt  algebras over a  field of prime characteristic \cite{Yao};
        \item[-] solvable Lie algebras of maximal rank \cite{KOK22};
        \item[-] Cayley algebras \cite{AEK1};
        \item[-] finite dimensional semi-simple Leibniz algebras over an algebraically closed field of characteristic zero \cite{KKY};
        \item[-] locally finite split simple Lie algebras over a field of characteristic zero \cite{AKYu1};
        \item[-] the Schr\"{o}dinger algebras \cite{AlaYus};
        \item[-] Lie superalgebra $q(n)$ \cite{ChWD};
        \item[-] conformal Galilei algebras \cite{AlaYus2}.
  	\end{itemize}	
	
	On the other hand, some algebras (in most cases close to nilpotent algebras) admit pure local derivations, that is, local derivations which  are not derivations. Below a short list of some classes of algebras which admit pure local derivations:
	\begin{itemize}
		\item[-] the algebra $\mathbb{C}(x)$ of rational functions \cite{Kad};
		\item[-] finite dimensional filiform Lie algebras \cite{AK};
        \item[-] $p$-filiform Leibniz algebras \cite{AKYu};
		\item[-] solvable Leibniz algebras with abelian nilradicals, which have a one dimensionial complementary space \cite{AKY20};
		\item[-] the  algebra of lower triangular $n\times n$-matrices \cite{E11};
        \item[-] the ternary Malcev algebra $M_8$ \cite{BIK};
        \item[-] direct sum null-filiform Leibniz algebras \cite{Adashev1}.
        \end{itemize}

     The Schr\"{o}dinger Lie group describes symmetries of the free particle Schr\"{o}dinger equation in \cite{Per}. The Lie algebra ${\rm {\mathbb S}}\left(n\right)$ in $\left(n+1\right)$-dimensional space-time of the Schr\"{o}dinger Lie group is called the Schr\"{o}dinger algebra, see \cite{DDM}. The Schr\"{o}dinger algebra is a non-semisimple Lie algebra and plays an important role in mathematical physics. The Lie algebra ${\rm {\mathbb S}}\left(1\right)$is one of the most essential case for $n=1$ and admits a universal $1$-dimensional central extension which is called the centrally extended Schr\"{o}dinger algebra or, simply, the Schr\"{o}dinger algebra, abusing the language.

The present paper is devoted to study local super-derivations on the super Schr\"{o}dinger algebra.

\section{Preliminaries}
In this section we recall definitions, symbols and establish some auxiliary results for
later use in this paper.

We begin by providing some conventions in this paper.
Denote the degree for $\mathbb{Z}_2$-graded
vector spaces $\mathcal{L} = \mathcal{L}_{\bar{0}}\oplus\mathcal{L}_{\bar{1}}$ by

 \begin{equation*}\label{Loc1} |x|=\left\{\begin{array}{ll}
 0,& \text{if}\ x\in\mathcal{L}_{\bar{0}}; \\[1mm]
 1,& \text{if}\ x\in\mathcal{L}_{\bar{1}}.\\[1mm]
   \end{array}\right.
 \end{equation*}

Elements in $\mathcal{L}_{\bar{0}}$ or $\mathcal{L}_{\bar{1}}$ are called homogeneous. We implicitely assume
that $x$ is as a $\mathbb{Z}_2$-homogeneous element when $|x|$ occurs in some expression. A linear map $g :\mathcal{L}\rightarrow\mathcal{L}$ is a homogeneous linear map of degree
$\alpha \in\mathbb{Z}_2$, i.e. $|g|=\alpha$, provided that

$$g(\mathcal{L}_{\beta})\subseteq\mathcal{L}_{\alpha+\beta},\ \forall \beta\in\mathbb{Z}_2.$$

Similarly, a bilinear map $\varphi : \mathcal{L}\times\mathcal{L}\rightarrow\mathcal{L}$ is a homogeneous bilinear map of degree $\alpha \in\mathbb{Z}_2$, i.e. $|\varphi|=\alpha$ if it satisfies that

$$\varphi(\mathcal{L}_{\beta},\mathcal{L}_{\gamma})\subseteq\mathcal{L}_{\alpha+\beta+\gamma},\ \forall\beta,\gamma\in\mathbb{Z}_2.$$

Moreover, if $|g|$ or $|\varphi|$ occurs in some expression, we regard it as a homogeneous map.

A Lie superalgebra (see \cite{Kac1977}) is a superspace $\mathcal{L} = \mathcal{L}_{\bar{0}} \oplus \mathcal{L}_{\bar{1}}$ with a bilinear mapping $ [., .] : \mathcal{L} \times \mathcal{L} \rightarrow \mathcal{L}$ satisfying the following identities:
\begin{enumerate}
\item $[\mathcal{L}_{\alpha}, \mathcal{L}_{\beta}] \subset \mathcal{L}_{\alpha+\beta}$, for $\alpha, \beta \in \mathbb{Z}_{2}$ ($\mathbb{Z}_{2}$-grading),
\item $[x, y] = -(-1)^{|x||y|} [y, x]$ (graded skew-symmetry),
\item $(-1)^{|x||z|} [x,[y, z]] + (-1)^{ |y| |x|} [y, [z, x]] + (-1)^{|z| |y|}[z,[ x, y]] = 0$ (graded Jacobi identity),
\end{enumerate}
for all $x, y, z \in \mathcal{L}$. It is clear that the  graded Jacobi identity  is equivalent to the equation
\begin{equation}\label{eq22}
[x,[y, z]]=  [[x,y], z] + (-1)^{|z| |x+y|}[[z,x],y].
\end{equation}

\begin{definition} \cite{CP} The $N=1$ super Schr\"{o}dinger algebra in $(1+1)$-dimensional spacetime $\mathcal{S}=\mathcal{S}_{\bar{0}}\oplus\mathcal{S}_{\bar{1}}$ is a 9-dimensional Lie superalgebra with the even part $\mathcal{S}_{\bar{0}} =span_{\mathbb{C}}\{e, f,h, p,q,z\}$ and the odd part $\mathcal{S}_{\bar{1}}=span_{\mathbb{C}}\{E,F,G\}$ satisfying the following non-vanishing Lie super brackets:
\begin{eqnarray*}
\begin{array}{lll}
[h,e]=2e,& [h,f]=-2f,&[e,f]=h,\\[1mm]
[h,p]=p,& [h,q]=-q, & [p,q]=z,\\[1mm]
[e,q]=p,& [p,f]=-q, & [h,E]=E,\\[1mm]
[h,F]=-F, & [e,F]=-E, & [f,E]=-F,\\[1mm]
[p,F]=G, & [q,E]=G,& [E,E]=2e,\\[1mm]
[F,F]=-2f,& [G,G]=z, & [E,F]=-h,\\[1mm]
[E,G]=-p, & [F,G]=q.\\[1mm]
\end{array}
\end{eqnarray*}
\end{definition}

\textbf{Remark 1.}  (1) The even part $\mathcal{S}_{\bar{0}} =span_{\mathbb{C}}\{e, f,h, p,q,z\}$ is the Schrodinger algebra in $(1+1)$-dimensional spacetime.

(2) It’s easy to see that $\mathcal{S}$ contains two subalgebras: the super Heisenberg subalgebra $H =
span_{\mathbb{C}}\{p,q,z,G\}$ and the Lie superalgebra $osp(1|2)$ spanned by $\{h,e, f,E,F\}.$ In particular, the
super Heisenberg subalgebra $H$ is also an ideal of $\mathcal{S}$. Then, $\mathcal{S}$ can be viewed as a semidirect
product $\mathcal{S}=H\rtimes osp(1|2).$

\begin{definition} Let $\mathcal{L}$ be a Lie superalgebra. A linear map $D :\mathcal{L}\rightarrow\mathcal{L}$ is called a super-derivation of $\mathcal{L}$ if
$$D([x,y])=[D(x),y]+(-1)^{|D||x|}[x,D(y)],\ \ \forall x,y\in\mathcal{L}.$$
\end{definition}

Write $Der_{\alpha}(\mathcal{L})$ for the set of all homogeneous super-derivations of degree $\alpha\in\mathbb{Z}_2$ of $\mathcal{L}$. Then
$Der(\mathcal{L})=Der_{\bar{0}}(\mathcal{L})\oplus Der_{\bar{1}}(\mathcal{L}).$ For arbitrary $x\in\mathcal{L}$, it is easy to see that the map $ad_x : \mathcal{L}\rightarrow\mathcal{L}$
defined by $ad_x(y)=[x,y]$ for all $y\in\mathcal{L}$ is a super-derivation of $\mathcal{L}$, which is termed an inner
super-derivation. The set of all inner super-derivations is denoted by $IDer(\mathcal{L}).$

\begin{definition}
A linear operator $\Delta$ is called a local super-derivation if for any $x \in \mathcal{L},$ there exists a super-derivation $D_x: \mathcal{L} \rightarrow \mathcal{L}$ (depending on $x$) such that
$\Delta(x) = D_x(x).$ The set of all local super-derivations on $\mathcal{L}$ we denote by $\mathrm{LocDer}(\mathcal{L}).$
\end{definition}

The following theorem gives a description of the super-derivation on the super Schr\"{o}dinger
algebra $\mathcal{S}.$

\begin{theorem}\label{Der1}\cite{YH}
$$\mathrm{Der}(\mathcal{S})=\mathrm{Ider}(\mathcal{S})\oplus\mathbb{C}\delta,$$
where $\delta$ are linear transformations of $\mathcal{S}$ and satisfy that
\begin{equation*}\begin{split}
\delta|_{osp(1|2)}=0 ,\ \delta(p)=p,\\ \delta(q)=q,\ \delta(z)=2z,\ \delta(G)=G.
\end{split}
\end{equation*}
\end{theorem}

\section{ Local super-derivations of the super Schr\"{o}dinger algebra $\mathcal{S}$}

In this section we investigate local super-derivations on the super Schr\"{o}dinger algebra $\mathcal{S}$.
The following theorem is the main result of this section.
\begin{theorem}\label{thm1}
Every local super-derivations on $\mathcal{S}$ is a super-derivation.
\end{theorem}

For the proof of this theorem we need several lemmas.

For any local super-derivations $\Delta$ on $\mathcal{S}$ we have by Theorem \ref{Der1}
$$\Delta(x)=[b_x,x]+\lambda_x\delta(x).$$
Put
$$b_x=b_{e,x}e+b_{f,x}f+b_{h,x}h+b_{p,x}p+b_{q,x}q+b_{z,x}z+b_{E,x}E+b_{F,x}F+b_{G,x}G\in\mathcal{S}.$$

Then we obtain
\begin{equation*}\begin{split}
\Delta(e)&=2b_{h,e}e-b_{f,e}h-b_{q,e}p+b_{F,e}E,\\
\Delta(f)&=-2b_{h,f}f+b_{e,f}h-b_{p,f}q+b_{E,f}F,\\
\Delta(h)&=-2b_{e,h}e+2b_{f,h}f-b_{p,h}p+b_{q,h}q-b_{E,h}E+b_{F,h}F,\\
\Delta(p)&=(b_{h,p}+\lambda_p)p+b_{f,p}q-b_{q,p}z-b_{F,p}G,\\
\Delta(q)&=b_{e,q}p+(-b_{h,q}+\lambda_q)q+b_{p,q}z-b_{E,q}G,\\
\Delta(z)&=2b_{z,z} z,\\
\Delta(E)&=2b_{E,E}e-b_{F,E}h+b_{G,E}p+b_{h,E}E-b_{f,E}F+b_{q,E}G,\\
\Delta(F)&=-2b_{F,F}f+b_{E,F}h-b_{G,F}q-b_{e,F}E-b_{h,F}F+b_{p,F}G,\\
\Delta(G)&=-b_{E,G}p+b_{F,G}q+b_{G,G}z+\lambda_G G.\\
\end{split}\end{equation*}

For the proof of this Theorem we need several Lemmata.

There exists an element $a=a_{e,x}e+a_{f,x}f+a_{h,x}h+a_{p,x}p+a_{q,x}q+a_{z,x}z+a_{E,x}E+a_{F,x}F+a_{G,x}G\in\mathcal{S}$ such that

$$\Delta(x)=[a,x]+\mu_x\delta(x).$$

In the following lemmas, we assume that $\Delta$ be a local super-derivation on $\mathcal{S}$ such that $\Delta(h+z)=0.$

\begin{lemma}\label{Lem1}
 $\Delta(h)=\Delta(z)=0.$
\end{lemma}
\begin{proof}
We considering,
$$0=\Delta(h+z)=\Delta(h)+\Delta(z)=-2b_{e,h}e+2b_{f,h}f-b_{p,h}p+b_{q,h}q-b_{E,h}E+b_{F,h}F+2b_{z,z}z,$$
which implies $b_{e,h}=b_{f,h}=b_{p,h}=b_{q,h}=b_{E,h}=b_{F,h}=b_{z,z}=0.$ Then $\Delta(h)=\Delta(z)=0.$
\end{proof}

\begin{lemma}
\begin{equation*}\begin{split}
\Delta(e)&=2b_{h,e}e-b_{q,e}p+b_{F,e}E,\\
\Delta(f)&=-2b_{h,f}f-b_{p,f}q+b_{E,f}F.\\
\end{split}\end{equation*}
\end{lemma}

\begin{proof}

We consider
\begin{equation*}\begin{split}
\Delta(h+e)&=[a,h+e]+\mu_{h+e}\delta(h+e)=\\
             &=[a_{e,h+e}e+a_{f,h+e}f+a_{h,h+e}h+a_{p,h+e}p+a_{q,h+e}q+a_{z,h+e}z+\\
             &+a_{E,h+e}E+a_{F,h+e}F+a_{G,h+e}G,h+e]+\mu_{h+e}\delta(h+e)=\\
             &=-2a_{e,h+e}e+2a_{f,h+e}f-a_{p,h+e}p+a_{q,h+e}q-a_{E,h+e}E+a_{F,h+e}F+\\
             &+2a_{h,h+e}e-a_{f,h+e}h-a_{q,h+e}p+a_{F,h+e}E.
\end{split}
\end{equation*}

On the other hand,
\begin{equation*}\begin{split}
\Delta(h+e)&=\Delta(h)+\Delta(e)=2b_{h,e}e-b_{f,e}h-b_{q,e}p+b_{F,e}E.
\end{split}
\end{equation*}

Comparing the coefficients at the basis elements $f$ and $h,$ we get $b_{f,e}=0.$

Similarly, considering

\begin{equation*}\begin{split}
\Delta(h+f)&=[a,h+f]+\mu_{h+f}\delta(h+f)=\\
             &=[a_{e,h+f}e+a_{f,h+f}f+a_{h,h+f}h+a_{p,h+f}p+a_{q,h+f}q+a_{z,h+f}z+\\
             &+a_{E,h+f}E+a_{F,h+f}F+a_{G,h+f}G,h+f]+\mu_{h+f}\delta(h+f)=\\
             &=-2a_{e,h+f}e+2a_{f,h+f}f-a_{p,h+f}p+a_{q,h+f}q-a_{E,h+f}E+a_{F,h+f}F\\
             &-2a_{h,h+f}f+a_{e,h+f}h-a_{p,h+f}q+a_{E,h+f}F.
\end{split}
\end{equation*}

On the other hand,
\begin{equation*}\begin{split}
\Delta(h+f)&=\Delta(h)+\Delta(f)=-2b_{h,f}f+b_{e,f}h-b_{p,f}q+b_{E,f}F.
\end{split}
\end{equation*}
Comparing the coefficients at the basis elements $e$ and $h,$ we get $b_{e,f}=0.$
\end{proof}

\begin{lemma}\label{lem4}

\begin{equation*}\begin{split}
\Delta(e)&=2b_{h,e}e-b_{q,e}p,\\
\Delta(f)&=-2b_{h,f}f-b_{p,f}q,\\
\Delta(p)&=(b_{h,p}+\lambda_p)p-b_{q,p}z,\\
\Delta(q)&=(-b_{h,q}+\lambda_q)q+b_{p,q}z.\\
\end{split}\end{equation*}
\end{lemma}

\begin{proof}
We consider
\begin{equation*}\begin{split}
\Delta(h+p)&=[a,h+p]+\mu_{h+p}\delta(h+p)=\\
             &=[a_{e,h+p}e+a_{f,h+p}f+a_{h,h+p}h+a_{p,h+p}p+a_{q,h+p}q+\\
             &+a_{z,h+p}z+a_{E,h+p}E+a_{F,h+p}F+a_{G,h+p}G,h+p]+
             \mu_{h+p}\delta(h+p)=\\
             &=-2a_{e,h+p}e+2a_{f,h+p}f-a_{p,h+p}p+a_{q,h+p}q-a_{E,h+p}E+a_{F,h+p}F+\\
             &+(a_{h,h+p}+\mu_{h+p})p+a_{f,h+p}q-a_{q,h+p}z-a_{F,h+p}G.
\end{split}
\end{equation*}

On the other hand,
\begin{equation*}\begin{split}
\Delta(h+p)&=\Delta(h)+\Delta(p)=
             (b_{h,p}+\lambda_p)p+b_{f,p}q-b_{q,p}z-b_{F,p}G.\\
\end{split}
\end{equation*}
Comparing the coefficients at the basis elements $F$ and $G,$ we get $b_{F,p}=0.$

We consider
\begin{equation*}\begin{split}
\Delta(h+q)&=[a,h+q]+\mu_{h+q}\delta(h+q)=\\
             &=[a_{e,h+q}e+a_{f,h+q}f+a_{h,h+q}h+a_{p,h+q}p+a_{q,h+q}q+\\
             &+a_{z,h+q}z+a_{E,h+q}E+a_{F,h+q}F+a_{G,h+q}G,h+q]+
             \mu_{h+p}\delta(h+p)=\\
             &=-2a_{e,h+q}e+2a_{f,h+q}f-a_{p,h+q}p+a_{q,h+q}q-a_{E,h+q}E+a_{F,h+q}F+\\
             &+a_{e,h+q}p+(-a_{h,h+q}+\mu_{h+q})q+a_{p,h+q}z-a_{E,h+q}G.
\end{split}
\end{equation*}

On the other hand,
\begin{equation*}\begin{split}
\Delta(h+q)&=\Delta(h)+\Delta(q)=
             b_{e,q}p+(-b_{h,q}+\lambda_q)q+b_{p,q}z-b_{E,q}G.\\
\end{split}
\end{equation*}
Comparing the coefficients at the basis elements $E$ and $G,$ we get $b_{E,q}=0.$

We consider
\begin{equation*}\begin{split}
\Delta(e+p)&=[a,e+p]+\mu_{e+p}\delta(e+p)=\\
             &=[a_{e,e+p}e+a_{f,e+p}f+a_{h,e+p}h+a_{p,e+p}p+a_{q,e+p}q+\\
             &+a_{z,e+p}z+a_{E,e+p}E+a_{F,e+p}F+a_{G,e+p}G,e+p]+
             \mu_{e+p}\delta(e+p)=\\
             &=2a_{h,e+p}e-a_{f,e+p}h-a_{q,e+p}p+a_{F,e+p}E+\\
             &+(a_{h,e+p}+\lambda_{e+p})p+a_{f,e+p}q-a_{q,e+p}z-a_{F,e+p}G.
\end{split}
\end{equation*}

On the other hand,
\begin{equation*}\begin{split}
\Delta(e+p)&=\Delta(e)+\Delta(p)=\\
             &=2b_{h,e}e-b_{q,e}p+b_{F,e}E+(b_{h,p}+\lambda_p)p+b_{f,p}q-b_{q,p}z.\\
\end{split}
\end{equation*}

Comparing the coefficients at the basis elements $h,\ q,\ E,$ and $G,$ we get $b_{F,e}=b_{f,p}=0.$

We consider
\begin{equation*}\begin{split}
\Delta(f+q)&=[a,f+q]+\mu_{f+q}\delta(f+q)=\\
             &=[a_{e,f+q}e+a_{f,f+q}f+a_{h,f+q}h+a_{p,f+q}p+a_{q,f+q}q+\\
             &+a_{z,f+q}z+a_{E,f+q}E+a_{F,f+q}F+a_{G,f+q}G,f+q]+
             \mu_{f+q}\delta(f+q)=\\
             &=-2a_{h,f+q}f+a_{e,f+q}h-a_{p,f+q}q+a_{E,f+q}F+\\
             &+a_{e,f+q}p+(-a_{h,f+q}+\lambda_{f+q})q+a_{p,f+q}z-a_{E,f+q}G.
\end{split}
\end{equation*}

On the other hand,
\begin{equation*}\begin{split}
\Delta(f+q)&=\Delta(f)+\Delta(q)=\\
             &=-2b_{h,f}f-b_{p,f}q+b_{E,f}F+b_{e,q}p+(-b_{h,q}+\lambda_q)q+b_{p,q}z.\\
\end{split}
\end{equation*}

Comparing the coefficients at the basis elements $h,\ p,\ F$ and $G,$ we get $b_{E,f}=b_{e,q}=0.$

\end{proof}

\begin{lemma}\label{lem4}
\begin{equation*}\begin{split}
\Delta(e)&=2b_{h,e}e-b_{q,e}p,\\
\Delta(f)&=-2b_{h,f}f-b_{q,e}q,\\
\Delta(E)&=2b_{E,E}e-b_{F,E}h+b_{G,E}p+b_{h,E}E-b_{f,E}F,\\
\Delta(F)&=-2b_{F,F}f+b_{E,F}h-b_{G,E}q-b_{e,F}E-b_{h,F}F,\\
\Delta(G)&=-b_{E,F}p+b_{F,E}q+b_{G,G}z+\lambda_G G.
\end{split}\end{equation*}
\end{lemma}

\begin{proof}

We consider
\begin{equation*}\begin{split}
\Delta(h+E)&=[a,h+E]+\mu_{h+E}\delta(h+E)=\\
             &=[a_{e,h+E}e+a_{f,h+E}f+a_{h,h+E}h+a_{p,h+E}p+a_{q,h+E}q+a_{z,h+E}z+\\
             &+a_{E,h+E}E+a_{F,h+E}F+a_{G,h+E}G,h+E]+\mu_{h+E}\delta(h+E)=\\
             &=-2a_{e,h+E}e+2a_{f,h+E}f-a_{p,h+E}p+a_{q,h+E}q-a_{E,h+E}E+a_{F,h+E}F+\\
             &+2a_{E,h+E}e-a_{F,h+E}h+a_{G,h+E}p+a_{h,h+E}E-a_{f,h+E}F+a_{q,h+E}G.
\end{split}
\end{equation*}

On the other hand,
\begin{equation*}\begin{split}
\Delta(h+E)&=\Delta(h)+\Delta(E)=2b_{E,E}e-b_{F,E}h+b_{G,E}p+b_{h,E}E-b_{f,E}F+b_{q,E}G.
\end{split}
\end{equation*}

Comparing the coefficients at the basis elements $q$ and $G,$ we get $b_{q,E}=0.$

Similarly, considering

\begin{equation*}\begin{split}
\Delta(h+F)&=[a,h+F]+\mu_{h+F}\delta(h+F)=\\
             &=[a_{e,h+F}e+a_{f,h+F}f+a_{h,h+F}h+a_{p,h+F}p+a_{q,h+F}q+a_{z,h+F}z+\\
             &+a_{E,h+F}E+a_{F,h+F}F+a_{G,h+F}G,h+F]+\mu_{h+F}\delta(h+F)=\\
             &=-2a_{e,h+F}e+2a_{f,h+F}f-a_{p,h+F}p+a_{q,h+F}q-a_{E,h+F}E+a_{F,h+F}F\\
             &-2a_{F,h+F}f+a_{E,h+F}h-a_{G,h+F}q-a_{e,h+F}E-a_{h,h+F}F+a_{p,h+F}G.
\end{split}
\end{equation*}

On the other hand,
\begin{equation*}\begin{split}
\Delta(h+F)&=\Delta(h)+\Delta(F)=-2b_{F,F}f+b_{E,F}h-b_{G,F}q-b_{e,F}E-b_{h,F}F+b_{p,F}G.
\end{split}
\end{equation*}
Comparing the coefficients at the basis elements $p$ and $G,$ we get $b_{p,F}=0.$

Considering

\begin{equation*}\begin{split}
\Delta(e-f+h)&=[a,e-f+h]+\mu_{e-f+h}\delta(e-f+h)=\\
             &=[a_{e,e-f+h}e+a_{f,e-f+h}f+a_{h,e-f+h}h+a_{p,e-f+h}p+a_{q,e-f+h}q+a_{z,e-f+h}z+\\
             &+a_{E,e-f+h}E+a_{F,e-f+h}F+a_{G,e-f+h}G,e-f+h]+\mu_{e-f+h}\delta(e-f+h)=\\
             &=2a_{h,e-f+h}e-a_{f,e-f+h}h-a_{q,e-f+h}p+a_{F,e-f+h}E+\\
             &+2a_{h,e-f+h}f-a_{e,e-f+h}h+a_{p,e-f+h}q-a_{E,e-f+h}F-\\
             &-2a_{e,e-f+h}e+2a_{f,e-f+h}f-a_{p,e-f+h}p+a_{q,e-f+h}q-a_{E,e-f+h}E+a_{F,e-f+h}F.\\.
\end{split}
\end{equation*}

On the other hand,
\begin{equation*}\begin{split}
\Delta(e-f+h)&=\Delta(e)-\Delta(f)=2b_{h,e}e-b_{q,e}p-2b_{h,f}f-b_{p,f}q.
\end{split}
\end{equation*}
Comparing the coefficients at the basis elements $p$ and $q,$ we get $b_{p,f}=b_{q,e}.$

We consider
\begin{equation*}\begin{split}
\Delta(E+F)&=[a,E+F]+\mu_{E+F}\delta(E+F)=\\
             &=[a_{e,E+F}e+a_{f,E+F}f+a_{h,E+F}h+a_{p,E+F}p+a_{q,E+F}q+a_{z,E+F}z+\\
             &+a_{E,E+F}E+a_{F,E+F}F+a_{G,E+F}G,E+F]+\mu_{E+F}\delta(E+F)=\\
             &=2a_{E,E+F}e-a_{F,E+F}h+a_{G,E+F}p+a_{h,E+F}E-a_{f,E+F}F+a_{q,E+F}G-\\
             &-2a_{F,E+F}f+a_{E,E+F}h-a_{G,E+F}q-a_{e,E+F}E-a_{h,E+F}F+a_{p,E+F}G.
\end{split}
\end{equation*}

On the other hand,
\begin{equation*}\begin{split}
\Delta(E+F)&=\Delta(E)+\Delta(F)=2b_{E,E}e-b_{F,E}h+b_{G,E}p+b_{h,E}E-b_{f,E}F-\\
          &-2b_{F,F}f+b_{E,F}h-b_{G,F}q-b_{e,F}E-b_{h,F}F.
\end{split}
\end{equation*}

Comparing the coefficients at the basis elements $p$ and $q,$ we get $b_{G,E}=b_{G,F}.$

We consider
\begin{equation*}\begin{split}
\Delta(E+G)&=[a,E+G]+\mu_{E+G}\delta(E+G)=\\
             &=[a_{e,E+G}e+a_{f,E+G}f+a_{h,E+G}h+a_{p,E+G}p+a_{q,E+G}q+a_{z,E+G}z+\\
             &+a_{E,E+G}E+a_{F,E+G}F+a_{G,E+G}G,E+G]+\mu_{E+G}\delta(E+G)=\\
             &=2a_{E,E+G}e-a_{F,E+G}h+a_{G,E+G}p+a_{h,E+G}E-a_{f,E+G}F+a_{q,E+G}G-\\
             &-a_{E,E+G}p+a_{F,E+G}q+a_{G,E+G}z+\mu_{E+G} G.
\end{split}
\end{equation*}

On the other hand,
\begin{equation*}\begin{split}
\Delta(E+G)&=\Delta(E)+\Delta(G)=2b_{E,E}e-b_{F,E}h+b_{G,E}p+b_{h,E}E-b_{f,E}F-\\
          &-b_{E,G}p+b_{F,G}q+b_{G,G}z+\lambda_G G.
\end{split}
\end{equation*}

Comparing the coefficients at the basis elements $h$ and $q,$ we get $b_{F,E}=b_{F,G}.$

We consider
\begin{equation*}\begin{split}
\Delta(F+G)&=[a,F+G]+\mu_{F+G}\delta(F+G)=\\
             &=[a_{e,F+G}e+a_{f,F+G}f+a_{h,F+G}h+a_{p,F+G}p+a_{q,F+G}q+a_{z,F+G}z+\\
             &+a_{E,F+G}E+a_{F,F+G}F+a_{G,F+G}G,F+G]+\mu_{F+G}\delta(F+G)=\\
             &=-2a_{F,F+G}f+a_{E,F+G}h-a_{G,F+G}q-a_{e,F+G}E-a_{h,F+G}F+a_{p,F+G}G-\\
             &-a_{E,F+G}p+a_{F,F+G}q+a_{G,F+G}z+\mu_{F+G} G.
\end{split}
\end{equation*}

On the other hand,
\begin{equation*}\begin{split}
\Delta(F+G)&=\Delta(F)+\Delta(G)=-2b_{F,F}f+b_{E,F}h-b_{G,F}q-b_{e,F}E-b_{h,F}F-\\
          &-b_{E,G}p+b_{F,G}q+b_{G,G}z+\lambda_G G.
\end{split}
\end{equation*}

Comparing the coefficients at the basis elements $h$ and $p,$ we get $b_{E,F}=b_{E,G}.$

\end{proof}

\begin{lemma}\label{lem5}
\begin{equation*}\begin{split}
\Delta(e)&=2b_{h,e}e,\\
\Delta(f)&=-2b_{h,e}f.\\
\end{split}
\end{equation*}
\end{lemma}

\begin{proof}

We consider take an element $y=e+f+E+F$

\begin{equation*}\begin{split}
\Delta(y)&=[a,y]+\mu_{y}\delta(y)=\\
             &=[a_{e,y}e+a_{f,y}f+a_{h,y}h+a_{p,y}p+a_{q,y}q+
             a_{z,y}z+a_{E,y}E+a_{F,y}F+a_{G,y}G,y]+\mu_{y}\delta(y)=\\
             &=2a_{h,y}e-a_{f,y}h-a_{q,y}p+a_{F,y}E-\\
             &-2a_{h,y}f+a_{e,y}h-a_{p,y}q+a_{E,y}F+\\
             &+2a_{E,y}e-a_{F,y}h+a_{G,y}p+a_{h,y}E-a_{f,y}F+a_{q,y}G-\\
             &-2a_{F,y}f+a_{E,y}h-a_{G,y}q-a_{e,y}E-a_{h,y}F+a_{p,y}G.\\
\end{split}
\end{equation*}

On the other hand,
\begin{equation*}\begin{split}
\Delta(y)&=\Delta(e)+\Delta(f)+\Delta(E)+\Delta(F)=2b_{h,e}e-b_{q,e}p-2b_{h,f}f-b_{q,e}q+\\
         &+2b_{E,E}e-b_{F,E}h+b_{G,E}p+b_{h,E}E-b_{f,E}F-\\
         &-2b_{F,F}f+b_{E,F}h-b_{G,E}q-b_{e,F}E-b_{h,F}F.
\end{split}
\end{equation*}

Comparing the coefficients at the basis elements $p,\ q$ and $G,$ we get

\begin{equation*}\label{yyyy}
  \begin{array}{ll}
    -b_{q,e}+b_{G,E}&=-a_{q,y}+a_{G,y},\\
    -b_{q,e}-b_{G,E}&=-a_{p,y}-a_{G,y},\\
    \ \ \ \ \ \ \ \ 0   &=a_{q,y}+a_{p,y},\\
    \end{array}
\end{equation*}
which implies $b_{q,e}=0.$

We consider

\begin{equation*}\begin{split}
\Delta(e+f)&=[a,e+f]+\mu_{e+f}\delta(e+f)=\\
             &=[a_{e,e+f}e+a_{f,e+f}f+a_{h,e+f}h+a_{p,e+f}p+a_{q,e+f}q+\\
             &+a_{z,e+f}z+a_{E,e+f}E+a_{F,e+f}F+a_{G,e+f}G,e+f]+\mu_{e+f}\delta(e+f)=\\
             &=2a_{h,e+f}e-a_{f,e+f}h-a_{q,e+f}p+a_{F,e+f}E-\\
             &-2a_{h,e+f}f+a_{e,e+f}h-a_{p,e+f}q+a_{E,e+f}F.
\end{split}
\end{equation*}

On the other hand,
\begin{equation*}\begin{split}
\Delta(e+f)&=\Delta(e)+\Delta(f)=2b_{h,e}e-2b_{h,f}f.\\
\end{split}
\end{equation*}

Comparing the coefficients at the basis elements $e$ and $f,$ we get $b_{h,e}=b_{h,f}.$

\end{proof}

\begin{lemma}\label{lem6}
\begin{equation*}\begin{split}
\Delta(p)&=(b_{h,p}+\lambda_p)p,\\
\Delta(q)&=(-b_{h,q}+\lambda_q)q,\\
\Delta(E)&=b_{G,E}p+b_{h,e}E,\\
\Delta(F)&=-b_{G,E}q-b_{h,e}F,\\
\Delta(G)&=b_{G,G}z+\lambda_G G.
\end{split}\end{equation*}
\end{lemma}
\begin{proof}
We consider

\begin{equation*}\begin{split}
\Delta(e+G)&=[a,e+G]+\mu_{e+G}\delta(e+G)=\\
             &=[a_{e,e+G}e+a_{f,e+G}f+a_{h,e+G}h+a_{p,e+G}p+a_{q,e+G}q+\\
             &+a_{z,e+G}z+a_{E,e+G}E+a_{F,e+G}F+a_{G,e+G}G,e+G]+\mu_{e+G}\delta(e+G)=\\
             &=2a_{h,e+G}e-a_{f,e+G}h-a_{q,e+G}p+a_{F,e+G}E-\\
             &-a_{E,e+G}p+a_{F,e+G}q+a_{G,e+G}z+\mu_{e+G}G.
\end{split}
\end{equation*}

On the other hand,
\begin{equation*}\begin{split}
\Delta(e+G)&=\Delta(e)+\Delta(G)=2b_{h,e}e-b_{E,F}p+b_{F,E}q+b_{G,G}z+\lambda_G G.\\
\end{split}
\end{equation*}

Comparing the coefficients at the basis elements $E$ and $q,$ we get $b_{F,E}=0.$

We consider

\begin{equation*}\begin{split}
\Delta(f+G)&=[a,f+G]+\mu_{f+G}\delta(f+G)=\\
             &=[a_{e,f+G}e+a_{f,f+G}f+a_{h,f+G}h+a_{p,f+G}p+a_{q,f+G}q+\\
             &+a_{z,f+G}z+a_{E,f+G}E+a_{F,f+G}F+a_{G,f+G}G,f+G]+\mu_{f+G}\delta(f+G)=\\
             &=-2a_{h,f+G}f+a_{e,f+G}h-a_{p,f+G}q+a_{E,f+G}F\\
             &-a_{E,f+G}p+a_{F,f+G}q+a_{G,f+G}z+\mu_{f+G}G.
\end{split}
\end{equation*}

On the other hand,
\begin{equation*}\begin{split}
\Delta(f+G)&=\Delta(f)+\Delta(G)=-2b_{h,f}f-b_{E,F}p+b_{G,G}z+\lambda_G G.\\
\end{split}
\end{equation*}

Comparing the coefficients at the basis elements $F$ and $p,$ we get $b_{E,F}=0.$

We consider take an element $y=E+F+G$
\begin{equation*}\begin{split}
\Delta(y)&=[a,y]+\mu_{y}\delta(y)=\\
             &=[a_{e,y}e+a_{f,y}f+a_{h,y}h+a_{p,y}p+a_{q,y}q+a_{z,y}z+a_{E,y}E+a_{F,y}F+a_{G,y}G,y]+\mu_{y}\delta(y)=\\
             &=2a_{E,y}e-a_{F,y}h+a_{G,y}p+a_{h,y}E-a_{f,y}F+a_{q,y}G-\\
             &-2a_{F,y}f+a_{E,y}h-a_{G,y}q-a_{e,y}E-a_{h,y}F+a_{p,y}G-\\
             &-a_{E,y}p+a_{F,y}q+a_{G,y}z+\mu_y G.
\end{split}
\end{equation*}

On the other hand,
\begin{equation*}\begin{split}
\Delta(y)&=\Delta(E)+\Delta(F)+\Delta(G)=2b_{E,E}e+b_{G,E}p+b_{h,E}E-b_{f,E}F-\\
         &-2b_{F,F}f-b_{G,E}q-b_{e,F}E-b_{h,F}F+b_{G,G}z+\lambda_G G.
\end{split}
\end{equation*}

Comparing the coefficients at the basis elements $e,\ h$ and $f,$ we get $b_{E,E}=b_{F,F}.$

We consider

\begin{equation*}\begin{split}
\Delta(e+p+E)&=[a,e+p+E]+\mu_{e+p+E}\delta(e+p+E)=\\
             &=[a_{e,e+p+E}e+a_{f,e+p+E}f+a_{h,e+p+E}h+a_{p,e+p+E}p+a_{q,e+p+E}q+\\
             &+a_{z,e+p+E}z+a_{E,e+p+E}E+a_{F,e+p+E}F+a_{G,e+p+E}G,e+p+E]+\\
             &+\mu_{e+p+E}\delta(e+p+E)=2a_{h,e+p+E}e-a_{f,e+p+E}h-\\
             &-a_{q,e+p+E}p+a_{F,e+p+E}E+(a_{h,e+p+E}+\mu_{e+p+E})p+\\
             &+a_{f,e+p+E}q-a_{q,e+p+E}z-a_{F,e+p+E}G+2a_{E,e+p+E}e-a_{F,e+p+E}h+\\
             &+a_{G,e+p+E}p+a_{h,e+p+E}E-a_{f,e+p+E}F+a_{q,e+p+E}G.
\end{split}
\end{equation*}

On the other hand,
\begin{equation*}\begin{split}
\Delta(e+p+E)&=\Delta(e)+\Delta(p)+\Delta(E)=2b_{h,e}e+\\
             &+(b_{h,p}+\lambda_p)p-b_{q,p}z+2b_{E,E}e+b_{G,E}p+b_{h,E}E-b_{f,E}F.\\
\end{split}
\end{equation*}

Comparing the coefficients at the basis elements $h,\ q,\ z,\ E,\ F$ and $G$ we get $b_{f,E}=b_{q,p}=0.$

We consider

\begin{equation*}\begin{split}
\Delta(q+F)&=[a,q+F]+\mu_{q+F}\delta(q+F)=\\
             &=[a_{e,q+F}e+a_{f,q+F}f+a_{h,q+F}h+a_{p,q+F}p+a_{q,q+F}q+\\
             &+a_{z,q+F}z+a_{E,q+F}E+a_{F,q+F}F+a_{G,q+F}G,q+F]+\mu_{q+F}\delta(q+F)=\\
             &=a_{e,q+F}p+(-a_{h,q+F}+\mu_{q+F})q+a_{p,q+F}z-a_{E,q+F}G-\\
             &-2a_{F,q+F}f+a_{E,q+F}h-a_{G,q+F}q-a_{e,q+F}E-a_{h,q+F}F+a_{p,q+F}G.
\end{split}
\end{equation*}

On the other hand,
\begin{equation*}\begin{split}
\Delta(q+F)&=\Delta(q)+\Delta(F)=(-b_{h,q}+\lambda_q)q+b_{p,q}z+\\
             &-2b_{F,F}f-b_{G,E}q-b_{e,F}E-b_{h,F}F.\\
\end{split}
\end{equation*}

Comparing the coefficients at the basis elements $h,\ p,\ z,\ E$ and $G,$ we get $b_{e,F}=b_{p,q}=0.$

We consider take an element $\nu=f+p+E$
\begin{equation*}\begin{split}
\Delta(\nu)&=[a,\nu]+\mu_{\nu}\delta(\nu)=\\
             &=[a_{e,\nu}e+a_{f,\nu}f+a_{h,\nu}h+a_{p,\nu}p+a_{q,\nu}q+a_{z,\nu}z+a_{E,\nu}E+a_{F,\nu}F+a_{G,\nu}G,\nu]+\mu_{\nu}\delta(\nu)=\\
             &=-2a_{h,\nu}f+a_{e,\nu}h-a_{p,\nu}q+a_{E,\nu}F+\\
             &+(a_{h,\nu}+\lambda_{\nu})p+a_{f,\nu}q-a_{q,\nu}z-a_{F,\nu}G+\\
             &+2a_{E,\nu}e-a_{F,\nu}h+a_{G,\nu}p+a_{h,\nu}E-a_{f,\nu}F+a_{q,\nu}G.
\end{split}
\end{equation*}

On the other hand,
\begin{equation*}\begin{split}
\Delta(\nu)&=\Delta(f)+\Delta(p)+\Delta(E)=-2b_{h,e}f+\\
         &+(b_{h,p}+\lambda_p)p-b_{q,p}z+2b_{E,E}e+b_{G,E}p+b_{h,E}E.
\end{split}
\end{equation*}

Comparing the coefficients at the basis elements $e,\ q$ and $F,$ we get $b_{E,E}=0.$

We consider

\begin{equation*}\begin{split}
\Delta(e+F)&=[a,e+F]+\mu_{e+F}\delta(e+F)=\\
             &=[a_{e,e+F}e+a_{f,e+F}f+a_{h,e+F}h+a_{p,e+F}p+a_{q,e+F}q+\\
             &+a_{z,e+F}z+a_{E,e+F}E+a_{F,e+F}F+a_{G,e+F},e+F]+\mu_{e+F}\delta(e+F)=\\
             &=2a_{h,e+F}e-a_{f,e+F}h-a_{q,e+F}p+a_{F,e+F}E-\\
             &-2a_{F,e+F}f+a_{E,e+F}h-a_{G,e+F}q-a_{e,e+F}E-a_{h,e+F}F+a_{p,e+F}G.
\end{split}
\end{equation*}

On the other hand,
\begin{equation*}\begin{split}
\Delta(e+F)&=\Delta(e)+\Delta(F)=2b_{h,e}e-b_{G,E}q-b_{h,F}F.\\
\end{split}
\end{equation*}

Comparing the coefficients at the basis elements $e,\ f,\ h,\ E$ and $F,$ we get $b_{h,F}=b_{h,e}.$

We consider

\begin{equation*}\begin{split}
\Delta(f+E)&=[a,f+E]+\mu_{f+E}\delta(f+E)=\\
             &=[a_{e,f+E}e+a_{f,f+E}f+a_{h,f+E}h+a_{p,f+E}p+a_{q,f+E}q+\\
             &+a_{z,f+E}z+a_{E,f+E}E+a_{F,f+E}F+a_{G,f+E},f+E]+\mu_{f+E}\delta(f+E)=\\
             &=-2a_{h,f+E}f+a_{e,f+E}h-a_{p,f+E}q+a_{E,f+E}F+\\
             &+2a_{E,f+E}e-a_{F,f+E}h+a_{G,f+E}p+a_{h,f+E}E-a_{f,f+E}F+a_{q,f+E}G.
\end{split}
\end{equation*}

On the other hand,
\begin{equation*}\begin{split}
\Delta(f+E)&=\Delta(f)+\Delta(E)=-2b_{h,e}f+b_{G,E}p+b_{h,E}E.\\
\end{split}
\end{equation*}

Comparing the coefficients at the basis elements $e,\ f,\ h,\ E$ and $F,$ we get $b_{h,E}=b_{h,e}.$

\end{proof}

\begin{lemma}\label{lem7}
\begin{equation*}\begin{split}
\Delta(p)&=b_{h,e}p,\\
\Delta(q)&=-b_{h,e}q,\\
\Delta(E)&=b_{G,E}p+b_{h,e}E,\\
\Delta(F)&=-b_{G,E}q-b_{h,e}F,\\
\Delta(G)&=b_{G,E}z.
\end{split}\end{equation*}
\end{lemma}
\begin{proof}

We consider take an element $y=E-F+G$
\begin{equation*}\begin{split}
\Delta(y)&=[a,y]+\mu_{y}\delta(y)=\\
             &=[a_{e,y}e+a_{f,y}f+a_{h,y}h+a_{p,y}p+a_{q,y}q+a_{z,y}z+a_{E,y}E+a_{F,y}F+a_{G,y},y]+\mu_{y}\delta(y)=\\
             &=2a_{E,y}e-a_{F,y}h+a_{G,y}p+a_{h,y}E-a_{f,y}F+a_{q,y}G+\\
             &+2a_{F,y}f-a_{E,y}h+a_{G,y}q+a_{e,y}E+a_{h,y}F-a_{p,y}G-\\
             &-a_{E,y}p+a_{F,y}q+a_{G,y}z+\mu_y G.
\end{split}
\end{equation*}

On the other hand,
\begin{equation*}\begin{split}
\Delta(y)&=\Delta(E)-\Delta(F)+\Delta(G)=b_{G,E}p+b_{h,e}E+\\
         &+b_{G,E}q+b_{h,e}F+b_{G,G}z+\lambda_G G.
\end{split}
\end{equation*}

Comparing the coefficients at the basis elements $e,\ p$ and $z,$ we get $b_{G,G}=b_{G,E}.$

We consider take an element $\nu=f+q-\frac{1}{2}z$
\begin{equation*}\begin{split}
\Delta(\nu)&=[a,\nu]+\mu_{\nu}\delta(\nu)=\\
             &=[a_{e,\nu}e+a_{f,\nu}f+a_{h,\nu}h+a_{p,\nu}p+a_{q,\nu}q+a_{z,\nu}z+a_{E,\nu}E+a_{F,\nu}F+a_{G,\nu}G,\nu]+
             \mu_{\nu}\delta(\nu)=\\
             &=-2a_{h,\nu}f+a_{e,\nu}h-a_{p,\nu}q+a_{E,\nu}F+\\
             &+a_{e,\nu}p+(-a_{h,\nu}+\mu_{\nu})q+a_{p,\nu}z-a_{E,\nu}G-\mu_{\nu} z.
\end{split}
\end{equation*}

On the other hand,
\begin{equation*}\begin{split}
\Delta(\nu)&=\Delta(f)+\Delta(p)-\frac{1}{2}\Delta(z)=-2b_{h,e}f+(-b_{h,q}+\lambda_q)q.
\end{split}
\end{equation*}

Comparing the coefficients at the basis elements $ f,\ q$ and $z,$ we get $-b_{h,q}+\lambda_q=-b_{h,e}.$

We consider take an element $\tau=e+p+\frac{1}{2}z$
\begin{equation*}\begin{split}
\Delta(\tau)&=[a,\tau]+\mu_{\tau}\delta(\tau)=\\
             &=[a_{e,\tau}e+a_{f,\tau}f+a_{h,\tau}h+a_{p,\tau}p+a_{q,\tau}q+a_{z,\tau}z+a_{E,\tau}E+a_{F,\tau}F+a_{G,\tau}G,\tau]+
             \mu_{\tau}\delta(\tau)=\\
             &=2a_{h,\tau}e-a_{f,\tau}h-a_{q,\tau}p+a_{F,\tau}E+\\
             &+(a_{h,\tau}+\mu_{\tau})p+a_{f,\tau}q-a_{q,\tau}z-a_{F,\tau}G+\mu_{\tau} z.
\end{split}
\end{equation*}

On the other hand,
\begin{equation*}\begin{split}
\Delta(\tau)&=\Delta(e)+\Delta(p)+\frac{1}{2}\Delta(z)=2b_{h,e}e+(b_{h,p}+\lambda_p)p.
\end{split}
\end{equation*}

Comparing the coefficients at the basis elements $ e,\ p$ and $z,$ we get $b_{h,p}+\lambda_p=b_{h,e}.$

We consider take an element $t=p+q+E-F+G$
\begin{equation*}\begin{split}
\Delta(t)&=[a,t]+\mu_{t}\delta(t)=\\
             &=[a_{e,t}e+a_{f,t}f+a_{h,t}h+a_{p,t}p+a_{q,t}q+a_{z,t}z+a_{E,t}E+a_{F,t}F+a_{G,t},t]+\mu_{t}\delta(t)=\\
             &=(a_{h,t}+\mu_{t})p+a_{f,t}q-a_{q,t}z-a_{F,t}G+\\
             &+a_{e,t}p+(-a_{h,t}+\mu_t)q+a_{p,t}z-a_{E,t}G+\\
             &+2a_{E,t}e-a_{F,t}h+a_{G,t}p+a_{h,t}E-a_{f,t}F+a_{q,t}G-\\
             &+2a_{F,t}f-a_{E,t}h+a_{G,t}q+a_{e,t}E+a_{h,t}F-a_{p,t}G-\\
             &-a_{E,t}p+a_{F,t}q+a_{G,t}z+\mu_t G.
\end{split}
\end{equation*}

On the other hand,
\begin{equation*}\begin{split}
\Delta(t)&=\Delta(p)+\Delta(q)+\Delta(E)-\Delta(F)+\Delta(G)=b_{h,e}p-\\
         &-b_{h,e}q+b_{G,E}p+b_{h,e}E-b_{G,E}q-b_{h,e}F+b_{G,G}z+\lambda_G G.
\end{split}
\end{equation*}

Comparing the coefficients at the basis elements $e,\ f,\  p,\ z,\ E$ and $G,$ we get $\lambda_G=0.$

\end{proof}

\begin{lemma}\label{Lem8}
$\Delta$ is an inner super-derivation on $\mathcal{S}$
\end{lemma}

\begin{proof} Let $x=x_{e}e+x_{f}f+x_{h}h+x_{p}p+x_{q}q+x_{z}z+x_{E}E+x_{F}F+x_{G}G\in\mathcal{S}.$

We consider, by Lemma \ref{Lem1}, Lemma \ref{lem5} and Lemma \ref{lem7}
\begin{equation*}\begin{split}\label{inn1}
\Delta(x)-&[b_{h,e}h+b_{G,E}G,x]=\Delta(x_{e}e+x_{f}f+x_{h}h+x_{p}p+x_{q}q+x_{z}z+x_{E}E+x_{F}F+x_{G}G)-\\
                      &-\left[b_{h,e}h+b_{G,E}G,x_{e}e+x_{f}f+x_{h}h+x_{p}p+x_{q}q+x_{z}z+x_{E}E+x_{F}F+x_{G}G\right]=\\
         &=(2x_{e}e-2x_{f}f+x_pp-x_qq+x_EE-x_FF)b_{h,e}+(x_Ep-x_Fq+x_Gz)b_{G,E}-\\
         &-(2x_{e}e-2x_{f}f+x_pp-x_qq+x_EE-x_FF)b_{h,e}-(x_Ep-x_Fq+x_Gz)b_{G,E}=0.
\end{split}
\end{equation*}
Then $\Delta$ is an inner super-derivation.

\end{proof}

Now we are in position to prove Theorem \ref{thm1}.

\textbf{Proof of Theorem \ref{thm1}}. Let $\nabla$ be a local super-derivation of $\mathcal{S}.$ For the element $h+z$ exists there a super-derivation
$D_{h+z}$ such that
$$\nabla(h+z)=D_{h+z}(h+z).$$

Set $\Delta=\nabla-D_{h+z}.$  Then $\Delta$ is a local super-derivation such that $\Delta(h+z)=0.$
By Lemma \ref{Lem8}, $\Delta(x) = ad_{b_{h,e}h+b_{G,E}G}(x).$
Thus $\nabla=D_{h+z}+ad_{b_{h,e}h+b_{G,E}G}.$  Then $\nabla$ is a super-derivation. The proof is complete.

\textbf{Declarations}

Conflict of interest The author declares no conflict of interest.

\end{document}